\documentclass{amsart}
\usepackage{graphics}
\usepackage{graphicx} 
\usepackage{psfrag} 
\usepackage{amssymb,a4wide,amsthm,amsmath,amsfonts}
\usepackage{mathrsfs}
\usepackage{color}

\newtheorem{theorem}{Theorem}[section]
\newtheorem{lemma}[theorem]{Lemma}

\theoremstyle{definition}
\newtheorem{definition}[theorem]{Definition}

\theoremstyle{remark}

\numberwithin{equation}{section}

\newcommand{\R}{\mathbb{R}}
\newcommand{\N}{\mathbb{N}}

\begin{document}
\title{Spaceability of sets of nowhere $L^q$ functions}

\author{ Pedro L. Kaufmann }
\address{Instituto de matem\'atica e estat\'istica, Universidade de S\~ao Paulo, Rua do Mat\~ao, 1010, CEP 05508-900, S\~ao Paulo,  Brazil}
\email{plkaufmann@gmail.com, leonardo@ime.usp.br}
\thanks{The first author was supported by CAPES, Research Grant PNPD 2256-2009.}

\author{ Leonardo Pellegrini}

\begin{abstract}

We say that a function $f:[0,1]\rightarrow \R$ is \emph{nowhere $L^q$} if, for each nonvoid open subset $U$ of $[0,1]$, the restriction $f|_U$ is not in $L^q(U)$. 
For a fixed $1\leq p <\infty$, we will show that the set 
$$
S_p\doteq \{f\in L^p[0,1]: f\mbox{ is nowhere $L^q$, for each }p<q\leq\infty\},  
$$
united with $\{0\}$, contains an isometric and complemented copy of $\ell_p$. In particular, this improves a result from \cite{bfps} (which in turn is an improvement of a result from \cite{mpps}), since $S_p$ turns out to be spaceable. In addition, our result is a generalization of one of the main results from \cite{gkp}. 

\end{abstract}
\maketitle

\section{Introduction}

This note is a contribution to the study of large linear structures within essentially non-linear sets of functions which satisfy certain special properties. Given a topological vector space $X$, it is by now standard to say that a subset $S\subset X$ is \emph{lineable} if $S\cup \{0\}$ contains an infinite dimensional subspace of $X$, and that $S$ is \emph{spaceable} if $S\cup \{0\}$ contains a \emph{closed} infinite dimensional subspace of $X$. 

Our main object of study will be the $L^p[0,1]$ spaces. For clearness, we will write $L^p$ instead of $L^p[0,1]$. Two results motivate the present work. Mu\~noz-Fern\'andez, Palmberg, Puglisi and Seoane-Sep�lveda showed the following: 

\begin{theorem}[from \cite{mpps}]

Let $1\leq p < q$. Then $L^p\setminus L^q$ is $\mathfrak{c}$-lineable. 
\label{motivo1}
\end{theorem}
Here, $\mathfrak{c}$ indicates that we can find a vector space with a Hamel basis of cardinality $\mathfrak{c}$ contained in $(L^p\setminus L^q)\cap\{0\}$. The other result is from Botelho, Diniz, F\'avaro and Pellegrino: 
\begin{theorem}[from \cite{bdfp}]

For every $p>0$, $\displaystyle{\ell_p \setminus\cup_{0<q<p} \ell_q}$ is spaceable. 
\label{motivo2}
\end{theorem}
In view of Theorems \ref{motivo1} and \ref{motivo2} it is natural to ask wether we can obtain a result for $L^p$ spaces analogue to Theorem \ref{motivo2};  Botelho, F\'avaro, Pellegrino and Seoane-Sep\'ulveda have obtained a positive answer to that question: 
\begin{theorem}[from \cite{bfps}]

For each $p>0$, $\displaystyle{L^p \setminus\cup_{q>p} L^q}$ is spaceable. 

\label{botelhoth}
\end{theorem}

In this work we will present an improvement of Theorem \ref{botelhoth}, one could say, to the \emph{local} level. We say that $f:[0,1]\rightarrow \R$ is \emph{nowhere} $L^q$, for some $1\leq q\leq \infty$, if for each nonvoid open subset $U$ of $[0,1]$, $f|_U$ is not in $L^q(U)$. 
For each $1\leq p <\infty$, we define 
\begin{eqnarray}
S_p\doteq \{f\in L^p([0,1]): f\mbox{ is nowhere $L^q$, for each }p<q\leq\infty\};
\label{mainset}   
\end{eqnarray}
then $S_p$ satisfies the following property, which is our main result: 
\begin{theorem}[Main]

For each $1\leq p <\infty$, $(S_p\cup\{0\})\subset L^p$ contains an isometric and complemented copy of $\ell_p$. In particular, $S_p$ is spaceable. 

\label{mainth}
\end{theorem}
Section \ref{secproof} will be dedicated to the proof of Theorem \ref{mainth}. Note that it only improves Theorem \ref{botelhoth} for $p\geq 1$. We decided to restrict ourselves to this case for the sake of uniformity of arguments; the careful reader can verify that the same construction can be used to show the spaceability of $S_p$ also for the quasi-Banach $p<1$ case. 
Note that the elements from $S_p$ are \emph{nowhere essentially bounded} (that is, nowhere $L^\infty$). It is worth mentioning that Theorem \ref{mainth} generalizes one of the main theorems of \cite{gkp}, which states the following:

\begin{theorem}[from \cite{gkp}]

The set $\mathcal{G}$ of $L^1$ functions which are nowhere essentially bounded is spaceable in $L^1$. 

\end{theorem}

In the same paper, it was shown that $\mathcal{G}$ is \emph{strongly $\mathfrak{c}$-algebrable}, that is, $\mathcal{G}\cup\{0\}$ contains a $\mathfrak{c}$-generated free algebra (with the pointwise multiplication). Let us mention that there are no nontrivial algebraic structures, not even 1-dimensional ones, within $S_p$, and even within $L^p\setminus L^q$, for $1\leq p <q < \infty$. In effect, if $f\in L^p\setminus L^q$ for such $p$ and $q$, then $f^n$ is not in $L^p$ for large enough $n$. 


\section{Proof of the main result}
\label{secproof}

We start by discriminating three objects that will be needed in order to construct a nice basic sequence of $L^p$ contained in $S_p$.  \\

\textbf{First: an almost disjoint family of  Cantor-built sets $(A_n)_n$.}
First, we recall the construction of the Cantor set of (Lebesgue) measure one half, which we will denote by $C$. We start with the closed unit interval (denote it by $C_1$), and remove the center open interval of lenght $1/4$, obtaining $C_2$, a disjoint union of two closed intervals. From each of them, we remove the center open interval of lenght $1/4^2$, obtaining $C_3$, the disjoint union of four closed intervals. Repeating the process and taking the intersection of all $C_n$, we obtain $C$. For our convenience, we will call a \emph{hole} of $C$ each open interval that was removed in each step of the construction of $C$. Let $t_n$ be the right endpoint of the interval which is the one more to the left, among the closed intervals that constitute $C_n$. The explicit formula for $t_n$ is $\frac{2+\sum_{j=2}^n 2^{j-2}}{2^{2n-1}}$, and for our purposes what is important about the sequence $t_n$ is the following:  
\begin{enumerate}
 \item $t_n < 2^{3-2n}$ and
\item $m(C\cap [0,t_n]) = 2^{n}$. 
\end{enumerate}
This information will be needed later on. We now proceed introducing an specific notation for sets which are essentially unions of sets similar to $C$.

\begin{definition}
\begin{enumerate}
 \item For each subinterval $I$ from $[0,1]$ and each subset  $A\subset [0,1]$, denote 
$$
A_I\doteq (b-a)A+a;
$$
\item a nonvoid subset $A\subset [0,1]$ of the form $\cup\{A_j:j\in\Gamma\}$ is said to be \emph{C-built} (with \emph{C-components $A_j$}) if there is an almost disjoint family $\{I_j:j\in\Gamma\}$ of subintervals of $[0,1]$ such that $A_j = C_{I_j}$, where $C$ is the Cantor set of measure $1/2$. 
\end{enumerate}
\label{defcbuilt}
\end{definition}  
According to the notation from \cite{gkp}, C-built sets are in particular \emph{Cantor-built}, which means, roughly speaking, that they are essentially unions of Cantor sets. The \emph{holes} of a set of the form $C_I$ are defined by the natural way: the set of all holes of $C_I$ is given by $\{J_I: J$ is a hole of $C\}$. With the notion of C-built sets in hands we can start partitioning $[0,1]$. First, put $A_1\doteq C$. $A_2$ is defined as follows:
$$
A_2\doteq \cup\{C_J: J\mbox{ is a hole of }A_1\}. 
$$
Note that $A_1$ and $A_2$ are almost disjoint, and that $m(A_2)=1/4$. Assuming that we have defined $A_n$, for some $n\geq 2$, then $A_{n+1}$ is defined inductively:
$$
A_{n+1}\doteq \cup\{C_J: J\mbox{ is a hole of a Cantor component of }A_n\}. 
$$
This way, we obtain a sequence $(A_n)_n$ of almost disjoint C-built sets satisfying 
\begin{enumerate}
\item $m(A_n)=2^{-n}$, and
\item for any given nonempty open subset $U$ of $[0,1]$, there exists an $n_0\in \N$ such that $A_{n}$ has a C-component contained in $U$, for each $n\geq n_0$. 
\end{enumerate}

\textbf{Second: a convenient element $h$ of $L^p\setminus \cup_{p<q\leq\infty} L^q$.} Let $(r_j)_j$ be a strictly decreasing sequence of real numbers converging to $p$. For each $j$, the function $h_j: x\mapsto x^{-1/r_j}$ is in $L^p$, but not in $L^q$ for any $r_j\leq q \leq \infty$. Define $\tilde{h}$ by 
\begin{eqnarray}
\tilde{h}\doteq \sum_{j=1}^\infty \frac{1}{2^j} \frac{h_j}{\|h_j\|_p}. 
\label{eqf}
\end{eqnarray}
$\tilde{h}$ is well-defined since the series converges absolutely and $L^p$ is complete. Then $\tilde{h}\in L^p\setminus \cup_{p<q\leq\infty} L^q$, and moreover: 

\begin{lemma}

$h\doteq\chi_C \tilde{h}\in L^p\setminus \cup_{p<q\leq\infty} L^q$.  

\end{lemma}

\textbf{Proof.} It is easily seen that $h\in L^p$. To show that, for a fixed  $p<q<\infty$ we have $h\not\in L^q$, consider  $j_0$ such that  $r_{j_0}< q$. Then
$$
\int |h|^q = \int_C |\tilde{h}|^q \geq 
\int_C \frac{1}{2^{j_0q}}\frac{h_{j_0}^q}{\|h_{j_0}\|_p^q} \geq 
\frac{1}{2^{j_0q}\|h_{j_0}\|_p^q } \int_C x^{-q/r_{j_0}}.
$$
We will show that the integral to the right converges to $+\infty$. In effect, denoting $s\doteq q/r_{j_0}>1$ and considering the sequence $t_n$ defined previously, we have that
\begin{eqnarray*}
\int_C x^{-q/r_{j_0}} \hspace{-0.2cm} &=&  \hspace{-0.2cm}\int_C x^{-s}\geq \int_{C\cap [0,t_n]} x^{-s}
\geq m(C\cap [0,t_n])\inf \{x^{-s}:x\in C\cap [0,t_n]\} \\
&\geq& \hspace{-0.2cm} m(C\cap [0,t_n])t_n^{-s} >  2^{-n} 2^{s(2n-3)} > 2^{(2s-1)n}. 
\end{eqnarray*}
Since the right handside tends to infinity along with $n$, it follows that $h\not\in L^q$. The remaining case, $q=\infty$, is automatically covered, since $L^\infty \subset L^q$, for each $q<\infty$. $\spadesuit$\\

\textbf{Third: a disjoint infinite family of strictly increasing sequences of positive integers, $\{(n_j^k)_j:k\in\N\}$.}  It is a simple task for the reader to convince herself/himself that such family exists. \\

Before proceeding with the proof of the main result, let us stablish the following notation: if  $I=[a,b]$ is a subinterval of $[0,1]$ and  $f:[0,1]\rightarrow \R$, we define $f_I:[0,1]\rightarrow \R$ by
$$
f_I(x) \doteq \left\{ \begin{array}{ll}
f(\frac{x-a}{b-a}), &\mbox{if } x\in I;\\
0, & \mbox{if } x\not\in I. 
\end{array}\right.
$$
Note that this notation is coherent with Definition \ref{defcbuilt} (1), since for any subset $B\subset [0,1]$, $(\chi_B)_I = \chi_{B_I}$. It is easily seen that, for any subinterval $I$ from $[0,1]$ and any $1\leq p\leq \infty$, $f$ is in $L^p$ if and only if $f_I$ is in $L^p$. In particular, for each subinterval $I$ from $[0,1]$, $h_I\in  L^p\setminus \cup_{p<q\leq\infty} L^q$.\\

\textbf{Proof (of Theorem \ref{mainth}).} Consider $(A_j)_j$, $h$ and $\{(n^k_j)_j:k\in\N\}$ as defined above. For a fixed $j\geq 2$, let $(B_l)_l$ be a sequence of all C-components of $A_j$, and consider their respective convex hulls $I_l$. Note that the $I_l$ are disjoint. Define 
$$
f_j \doteq \sum_{l=1}^\infty \frac{1}{2^l}\frac{h_{I_l}}{\|h_{I_l}\|_p}. 
$$
For $j=1$ just put $f_1\doteq \frac{h}{\|h\|_p}$. Then each $f_j$ is in $L^p$ with $\|f_j\|_p=1$, and $f_j$ is not in $L^q(U)$, for any $p<q\leq\infty$ and any open set $U$ containing a C-component of $A_j$. Note also, by the construction of $h$, that $f_j$ is zero outside of $A_j$. For each $k$, consider
$$
g_k \doteq \sum_{j=1}^\infty \frac{1}{2^j} f_{n_j^k}.
$$
Then $g_k \in S_p$ and $\|g_k\|_p=1$, and since the functions $g_k$ have almost
disjoint supports, it follows that $(g_k)_k$ is a complemented basic sequence in $L^p$, isometrically equivalent to the canonical basis of $\ell_p$. In particular, $\overline{span(\{g_k:k\in \N\})}$ is complemented in $L^p$ and isometrically isomorphic to $\ell_p$. 

To complete our proof, it remains to show that $\overline{span(\{g_k:k\in \N\})}\subset S_p\cup \{0\}$. Any given nonzero element of $\overline{span(\{g_k:k\in \N\})}$ is of the form 
$$
g=\sum_{k=1}^\infty a_k g_k,  
$$
where some $a_{k_0}$ is assumed to be nonzero. Since the functions $g_k$ have almost disjoint supports and $g_{k_0}$ is nowhere $L^q$ for any given $p<q\leq\infty$, it follows that $g\in S_p$. $\spadesuit$


\end{document}